\documentclass[11pt]{amsart}    

\usepackage{graphicx}
\usepackage{color}
 
\usepackage{amsfonts,amsmath,latexsym,amssymb,amsthm}
\usepackage{hyperref}
\usepackage{geometry}
\usepackage{txfonts}

\newtheorem{theorem}{Theorem}

\theoremstyle{remark}
\newtheorem{remark}{Remark}

\newtheorem{lemma}{Lemma}
 
\begin{document}

\title{Magnetic Schr\"{o}dinger operator with the potential supported in a curved two-dimensional strip}

\author{Juan Bory Reyes}
  
\address{ESIME-Zacatenco, Instituto Polit\'ecnico Nacional, M\'exico,\\ CDMX. 07738. M\'exico.} 
\email{juanboryreyes@yahoo.com}
 
\author{Baruch Schneider}
 
 \address{ Department of Mathematics, Faculty of Science, University of Ostrava \\ 30.~dubna 22, 70103 Ostrava, Czech Republic}
\email{baruch.schneider@osu.cz}
  
\author{Diana Schneiderov\'{a} (Barseghyan)}
\address{ Department of Mathematics, Faculty of Science, University of Ostrava\\ 30.~dubna 22, 70103 Ostrava, Czech Republic}
\email{diana.schneiderova@osu.cz}


\maketitle

\begin{abstract}
We consider the magnetic Schr\"odinger operator $H=(i \nabla +A)^2- V$ with a non-negative potential $V$ supported over a strip which is a local deformation of a straight one, and the magnetic field $B:=\mathrm{rot}(A)$ is assumed to be nonzero and local. We show that the magnetic field does not change the essential spectrum of this system, and investigate a sufficient condition for the discrete spectrum of $H$ to be empty.
\end{abstract}


\section{Introduction}
\label{s:intro}

It is well known that a magnetic field, even a local one, can significantly affect the behavior of waveguide systems, in particular the existence of a geometrically induced discrete spectrum. While a particle confined in a 
fixed-profile tube with a hard-wall boundary can exist in localized states whenever the tube is bent or locally deformed (and asymptotically straight), cf.~\cite{EK15} for a comprehensive review of quantum waveguide theory, the presence of a local magnetic field can destroy such a discrete spectrum. This is a consequence of a Hardy-type inequality proved by Ekholm and Kova\v{r}\'{\i}k \cite{EK05}.

In this work we consider two-dimensional magnetic
Schr\"{o}dinger operators with an attractive potential supported over a strip which is a local perturbation of a straight tube.   
The similar model but with zero magnetic field was considered in \cite{E20}. If a corresponding strip is not straight, but is straight outside a bounded region, or it is at least asymptotically straight in a suitable sense, then for potentials with some suitable properties the Hamiltonian has a nontrivial geometrically induced discrete spectrum below the threshold of the essential spectrum, which is stable with respect to the local perturbations of a strip. 
 
The analogous effect of the existence of bound states resulting from the geometry was observed in \cite{ES89}, \cite{EK051}, \cite{EKK08}.

Our topic of interest concerns the influence of external fields in a model described in \cite{E20}. It is well known that when the magnetic field is homogeneous, both the essential and the discrete spectrum change.
While a local field is a much weaker perturbation that does not change the essential spectrum, the question is what must happen for the discrete spectrum to be destroyed.

The similar effect of the magnetic field has been observed in a class of singular Schr\"{o}dinger operators, usually called leaky quantum wires, with attractive contact interaction supported by a curve in \cite{BE21}.

Let us define the operator considered in this manuscript. Let $\Gamma$ be an infinite and smooth planar curve without self-intersections, i.e. the graph of any smooth function. We define the strip $\Omega^a$ in the plane, of half-width $a > 0$, built over $\Gamma$ as $$\Omega^a:= \{z \in\mathbb{R}^2: \mathrm{dist}(z, \Gamma) < a\}.$$

The main object of our interest are magnetic Schr\"{o}dinger operators with a compactly supported magnetic field and an attractive potential $V$ over the strip $\Omega^a$. We introduce $V$ in terms of the one-dimensional function as follows: we consider a non-zero $\tilde{V}\ge0$ of $L^\infty(\mathbb{R})$ with $\mathrm{supp} \tilde{V}\subset [-a, a]$
and define 
$V : \Omega^a\to\mathbb{R}_+,\, V(x,y) = \tilde{V}(u)$, where $u=u(x,y)$ is the distance from point $(x,y)$ to $\Gamma$.

We will study the operator
\begin{equation}\label{Hamiltonian}H=(i \nabla +A)^2- V\,
\end{equation} 
where $A$ is a vector potential corresponding to the magnetic field $B$.

There are several ways to properly define the operator , the simplest one is to identify it with the
unique self-adjoint operator associated with the quadratic form
\begin{equation}\label{q}
q(f)=\int_{\mathbb{R}^2}|i\nabla f+ A f|^2\,d x\,d y- \int_{\mathbb{R}^2}V |f|^2\,d x\,d y,\quad f\in \mathcal{H}^1(\mathbb{R}^2).
\end{equation}

The following assumptions will be required throughout the paper

\begin{itemize} 
\item
the magnetic field $B = \mathrm{rot}\, A\in L_{\mathrm{loc}}^2(\mathbb{R}^2)$ is compactly supported with $\mathrm{supp} (B) \subset\mathcal{B}(0, c)$, where $\mathcal{B}(0, c)$ is the ball centered at the origin, with radius $c$ for some $c > 0$;

\item
outside of $\mathcal{B}(0, c)$ the curve $\Gamma$ coincides with the straight line $\Gamma_k := \{(x, k x)\}_{x\in\mathbb{R}},\,k\in \mathbb{R}$;

\item
we make no additional assumptions except for smoothness on the curve $\Gamma$ in $\mathcal{B}(0, c)$.
\end{itemize}

For the further description it is necessary to introduce the one-dimensional auxiliary operator acting in $L^2(\mathbb{R})$,
$$h= -\frac{d^2}{d x^2}- \tilde{V}(x),$$
with domain $\mathcal{H}^2(\mathbb{R})$.

In \cite{E20} it was shown that if $\Omega^a$ is a straight tube and the magnetic field is absent, then due to
the separation of variables, the spectrum of operator $H$ is purely essential and coincides with $[e, \infty)$, where  
$e:=\mathrm{inf} \,\sigma (h)$. In the same work it was shown that the local perturbation of the straight tube does not change the essential spectrum. However, the discrete spectrum of such an operator is nonempty if the potential is large and the strip is narrow enough. We will show that, under some suitable assumptions, the presence of a local magnetic field can destroy the discrete spectrum. 

In Section 2 we present our main results. In fact, we prove the stability of the
essential spectrum of the operator (\ref{Hamiltonian}) with respect to the situation when the strip is a local deformation of a straight tube and the magnetic field is non-zero and local, and establish a
a sufficient condition for the absence of its discrete spectrum. As an application we give the example of paper \cite{E20}
where the curve deformation produces the discrete spectrum and we show that the local magnetic field destroys the eigenvalues.   

\section{Main results}\label{main}

The following notations will be used. Let $\mathcal{B}(0, r)$ be the open ball centered on the origin with radius $r$, and 
\begin{equation}\label{flux}
\Phi(r)= \frac{1}{2\pi}\int_{\mathcal{B}(0, r)}B(x, y)\,d x\,d y
\end{equation}
is the flux of $B$ through the ball $\mathcal{B}(0, r)$.

For any interval $I\subset\mathbb{R}$ let $d(I)$ be the length of $I$.

First we consider the situation when the curve $\Gamma$ coincides with the coordinate axis $\Gamma_0=\{(x, 0)\}_{x\in \mathbb{R}}$ outside of the ball $\mathcal{B}(0, c)$. The following theorems present our first and second results: 

\begin{theorem}\label{th.11.}
Suppose that, the curve $\Gamma$ coincides with the coordinate axis $\Gamma_0$ outside $\mathcal{B}(0, c),\, c > 0$. Under the assumptions made in the introduction, the essential spectrum of the operator $H$ coincides with the half-line $\big[e, \infty\big)$.
\end{theorem}

\begin{theorem}\label{th.22}
Let the assumptions of Theorem \ref{th.11.} hold, and further  suppose that the function $\Phi(r)$ is not identically zero for $r\in (0, c)$. Then there exists a constant $\mathcal{C}  > 0$ depending on the magnetic field such that the discrete spectrum of the operator $H$ is empty, provided that
\begin{equation}\label{assumption1}\left\|\tilde{V}(u)- \tilde{V}(y)\right\|_{L^\infty(\Omega^a)} \le \mathcal{C} \left(\delta(\Omega^a)\right)^{-1}\,,\end{equation}
where for any point $(x, y)\in\Omega^a$ we denote by $u=u(x, y)$ the distance from $(x, y)$ to $\Gamma$, and
$$
\delta(\Omega^a)=\underset{x_0\in(-c, c)}{\mathrm{max}} d\left(\{x=x_0\}\cap\Omega^a\right)\,.
$$
\end{theorem}

Before stating our results for a more general situation, we need to make the following two remarks:  

\begin{remark}\label{remark}

Suppose the curve $\Gamma$ coincides with $\Gamma_k$ outside the ball $\mathcal{B}(0, c)$. 
Let $\alpha_k\ne 0$ be the angle between $\Gamma_k$ and the coordinate axis $\Gamma_0=\{(x, 0)\}_{x\in\mathbb{R}}$ and let $\mathcal{U}$ be the unitary transformation on $\mathbb{R}^2$ to itself defined by
$$\mathcal{U}(x, y):=(x\cos\alpha_k- y\sin\alpha_k, x\sin\alpha_k+y\cos\alpha_k).$$

Now we define the unitary transformation on $L^2(\mathbb{R}^2)$ to itself as follows 
$$
(Uf)(x, y):= f(\mathcal{U}(x, y)) =f(x\cos\alpha_k- y\sin\alpha_k, x\sin\alpha_k+y\cos\alpha_k)\,.
$$

It is easy to check that for any $f\in \mathcal{H}^1(\mathbb{R}^2)$ the quadratic form $q$ can be rewritten as follows

\begin{eqnarray*}
q(f)=\int_{\mathbb{R}^2}\biggl(\left|\frac{\partial (Uf)}{\partial x}\right|^2+ i\mathcal{A}\left(-\frac{\partial (Uf)}{\partial x}\overline{(Uf)}+ \frac{\partial \overline{(Uf)}}{\partial x} (Uf)\right)\cos\alpha_k\\
+\left|\frac{\partial (Uf)}{\partial y}\right|^2+  i\mathcal{A}\left(\frac{\partial (Uf)}{\partial y}\overline{(Uf)}- \frac{\partial \overline{(Uf)}}{\partial y} (Uf)\right)\sin\alpha_k+ \mathcal{A}^2|Uf|^2\biggr)\,d xd y- \int_{\mathbb{R}^2} \mathcal{V}|Uf|^2\,d x\,d y\,,
\end{eqnarray*}
where
\begin{eqnarray*}\nonumber&\mathcal{A}(x, y)&\\&=\left(-\left(\int_0^{x\sin\alpha_k+y\cos\alpha_k}B(x\cos\alpha_k- y\sin\alpha_k, t)\,d t\right) \cos\alpha_k, \,\left(\int_0^{x\sin\alpha_k+y\cos\alpha_k}B(x\cos\alpha_k- y\sin\alpha_k, t)\,d t\right)\sin\alpha_k\right)\,,&\\\label{V}&
\text{and}\, \,\mathcal{V}(x, y)= V(x\cos\alpha_k- y\sin\alpha_k, x\sin\alpha_k+y\cos\alpha_k).&
\end{eqnarray*}

Then the right-hand side of the above expression equals  
$$q_0(Uf)=\int_{\mathbb{R}^2}|i\nabla(Uf)+ \mathcal{A}(Uf)|^2\,d x\,d y- \int_{\mathbb{R}^2} \mathcal{V}|Uf|^2\,d x\,d y,
$$
where $q_0$ is the quadratic form of the operator $H_0= (i \nabla+ \mathcal{A})^2- \mathcal{V}$ defined on its domain in $L^2(\mathbb{R}^2)$.

Combining this with the fact that the domains of the forms $q$ and $q_0$ coincide, we get the unitary
\begin{equation}\label{equivalence}
H=UH_0U^{-1}\,.
\end{equation} 

The various properties take place (we omit their proofs here for the sake of simplicity).
\begin{itemize}

\item 
Let $\Omega^a$ be the strip built over $\Gamma$ as usual.  Then $\mathcal{U}^{-1}\Omega^a$ coincides with the straight tube $\mathbb{R}\times(-a, a)$ outside of the ball $\mathcal{B}(0, c)$. Furthermore, $\mathcal{U}^{-1}\Omega^a$ coincides with the support of the potential $\mathcal{V}$.

\item The magnetic potential $\mathcal{A}$ is satisfied by   
\begin{eqnarray*}
\mathrm{rot}\,(\mathcal{A})= \frac{\partial}{\partial x}\left(\int_0^{x\sin\alpha_k+y\cos\alpha_k}B(x\cos\alpha_k-y\sin\alpha_k, t)\,d t\right)\sin\alpha_k\\+  \frac{\partial}{\partial y}\left(\int_0^{x\sin\alpha_k+y\cos\alpha_k}B(x\cos\alpha_k-y\sin\alpha_k, t)\,d t\right)\cos\alpha_k\\=
B(x\cos\alpha_k- y\sin\alpha_k, x\sin\alpha_k+y\cos\alpha_k)\,.
\end{eqnarray*}

\item
Due to the above statement the magnetic fluxes corresponding to the magnetic potentials $\mathcal{A}$ and $A$ are identical.

\item
Since the distance $u=u(x, y)$ from any point 
$(x, y)\in\Omega^a$ to $\Gamma$ is equal to the distance from point $\mathcal{U}^{-1}(x, y)\in\mathcal{U}^{-1}\Omega^a$ to $\mathcal{U}^{-1}\Gamma$, then for any $(x_1, y_1)\in\mathcal{U}^{-1}\Omega^a$ we have
 \begin{eqnarray*}\mathcal{V}(x_1, y_1)= V(x_1\cos\alpha_k- y_1\sin\alpha_k, x_1\sin\alpha_k+y_1\cos\alpha_k)\\
 =\tilde{V}(u(x_1\cos\alpha_k- y_1\sin\alpha_k, x_1\sin\alpha_k+y_1\cos\alpha_k))= \tilde{V}(w(x_1, y_1))\,,\end{eqnarray*}
where $w=w(x_1, y_1)$ is the distance from $(x_1, y_1)$ to $\mathcal{U}^{-1}\Gamma$.

\item
The following expression
$$
\left\|\tilde{V}(u)- \tilde{V}(d_{x, y})\right\|_{L^\infty(\Omega^a)}= \left\|\tilde{V}(u(x, y))- \tilde{V}(d_{x, y})\right\|_{L^\infty(\Omega^a)}
$$ 
where $d_{(x, y)}$ is the distance from point $(x, y)\in\Omega^a$ to line $\Gamma_k$, coincides with 
$$
\left\|\tilde{V}(w)- \tilde{V}(y_1)\right\|_{L^\infty(U^{-1}\Omega^a)}\,,
$$
where $w=w(x_1, y_1)$ is the distance from $(x_1, y_1)$ to $\mathcal{U}^{-1}\Gamma$.

\item

Let 
$$
\delta_k(\Omega^a):=\underset{x_0\in(-c, c)}{\mathrm{max}}d(\gamma_{x_0}\cap\Omega^a),
$$
where $\gamma_{x_0}$ is the line orthogonal to $\Gamma_k$ which heats $\Gamma_k$ at the point $(x_0, k x_0)$. Then
$$
\delta_k(\Omega^a)= \delta(\mathcal{U}^{-1}\Omega^a)$$ with $\delta(\cdot)$ is defined in Theorem(\ref{th.22}). 
\end{itemize}

\end{remark}

In view of the unitary equivalence (\ref{equivalence}) and the other statements of Remark \ref{remark}, we are able to extend the Theorems \ref{th.11.}-\ref{th.22} in the case when the curve $\Gamma$ coincides with $\Gamma_k$ outside the ball $\mathcal{B}(0, c)$. Below we present our second and third main results, using the same notation for the operator $H$ as in the Theorems \ref{th.11.}-\ref{th.22}.  

\begin{theorem}\label{th.1}
Suppose that, the curve $\Gamma$ coincides with the line $\Gamma_k,\,k\in\mathbb{R}$ outside $\mathcal{B}(0, c),\, c > 0$. Then under the assumptions of Theorem \ref{th.11.}; the essential spectrum of the operator $H$ coincides with the half-line $\big[e, \infty\big)$.
\end{theorem}
 
\begin{theorem}\label{th.2}
Let the assumptions of Theorem \ref{th.1} hold, and further suppose that the function $\Phi(r)$ is not identically zero for $r\in (0, c)$. Then there exists a constant $\mathcal{C}  > 0$ depending on the magnetic field such that the discrete spectrum of the operator $H$ is empty, provided that
\begin{equation}\label{assumption}\left\|\tilde{V}(u)- \tilde{V}\left(d_{(x, y)}\right)\right\|_{L^\infty(\Omega^a)} \le  \mathcal{C} \left(\delta_k(\Omega^a)\right)^{-1}\,,\end{equation}
where $u=u(x, y)$,  $d_{(x, y)},$ and $\delta_k(\cdot)$ are defined in Remark \ref{remark}.\end{theorem}

\section{Proofs}
\label{proofs}

Let us begin with the \emph{proof of Theorem~\ref{th.11.}}. First we will show that the essential spectrum of $H$ contains the half-line $\big[e, \infty\big)$. The proof essentially uses the Weyl criterion \cite[Thm~VII.12]{RS81}. Put $\lambda=e+p^2,\,p\in\mathbb{R}$, and consider the sequence of vectors
$$
\psi_k:\: \psi_k(x, y)=\frac{1}{\sqrt{k}}\,f(y)\,\mathrm{e}^{i p x}\,\chi\left(\frac{x}{k}\right),
$$
where $f$ is the ground state eigenfunction of $h$, $\chi$ is a smooth function with support in $(1, 2)$, and $k\in \mathbb{N}$.

It is easy to check that $\psi_k\in \mathrm{Dom}(H)$; the idea is to use the fact that there is a part of the plane where the magnetic field has no influence. In fact, by choosing the Landau gauge for the vector potential $A$,i.e. by setting 
$A(x, y)= \Big(-\int_0^y 
B(x, t)\,\mathrm{d}t,0\Big)$, we get $A(x, y)=0$ for $|x|>c$ thanks to the fact that $\mathrm{supp}\,(B) \subset\mathcal{B}(0, c)$. Then, for sufficiently large values of $k$
\begin{equation}\label{k}\int_{\mathbb{R}^2}|H \psi_k- \lambda \psi_k|^2\,d x\,d y= \int_{\mathbb{R}^2}|-\Delta \psi_k- V \psi_k- \lambda \psi_k|^2\,d x\,d y.\end{equation}

One has
\begin{align*}&
\frac{\partial^2\psi_k}{\partial y^2}=\frac{1}{\sqrt{k}}\,f''(y)\,\mathrm{e}^{i p x}\,\chi\left(\frac{x}{k}\right)\,,\\&
\frac{\partial^2\psi_k}{\partial x^2}=\frac{1}{\sqrt{k}}\,\left(-p^2 f(y)\,\chi\left(\frac{x}{k}\right)+ \frac{2 i p}{k} f(y)\,\chi'\left(\frac{x}{k}\right)+ \frac{1}{k^2}
f(y)\,\chi''\left(\frac{x}{k}\right)\right) \,\mathrm{e}^{i p x}\,.
\end{align*}

Hence in view of (\ref{k}) and the fact that $V(x, y)= \tilde{V}(y)$ for $|x|>c$, we have
\begin{gather*}
\int_{\mathbb{R}^2}\left|H \psi_k- \lambda \psi_k\right|^2\,d x\,d y\\\nonumber= \frac{1}{k} \int_k^{2k}\int_{\mathbb{R}}\left|-\frac{1}{k^2} f(y)\,\chi''\left(\frac{x}{k}\right)- f''(y)\,\chi\left(\frac{x}{k}\right)- \frac{2 i p}{k} f(y)\,\chi'\left(\frac{x}{k}\right)-
\tilde{V}(y)f(y)\,\chi\left(\frac{x}{k}\right)- e f(y)\,\chi\left(\frac{x}{k}\right)\right|^2\,d x\,d y\,.\end{gather*}

Since $f$ is the ground state eigenfunction of $h$ corresponding to eigenvalue $e$ the above expression implies
\begin{align}\nonumber&
\int_{\mathbb{R}^2}\left|H \psi_k- \lambda \psi_k\right|^2\,dx\,dy\\\nonumber&= \frac{1}{k} \int_k^{2k}\int_{\mathbb{R}}\biggr|\frac{2 i p}{k} f(y)\,\chi'\left(\frac{x}{k}\right)+ \frac{1}{k^2} f(y)\,\chi''\left(\frac{x}{k}\right)\biggr|^2\,d x\,d y
&\\\nonumber\le &\frac{4 p^2}{k^3} \int_k^{2k}\int_{\mathbb{R}} f(y)^2\,\left(\chi'\left(\frac{x}{k}\right)\right)^2\,d x\,d y+ \frac{2}{k^5}  \int_k^{2k}\int_{\mathbb{R}} f(y)^2\,\left(\chi''\left(\frac{x}{k}\right)\right)^2\,d x\,d y&\\
\nonumber&= \frac{4 p^2}{k^2} \int_{\mathbb{R}}f(y)^2\,d y \int_1^2 (\chi'(t))^2\,d x\,d y+ \frac{2}{k^4} \int_{\mathbb{R}} f(y)^2\,dy  \int_1^2(\chi''(t))^2\,d t
=\mathcal{O}\left(\frac{1}{k^2}\right),&\end{align}
and since  $\|\psi_k\|^2$ are independent of $k$, we infer that $\lambda =e+p^2\in \sigma(H)$. Moreover, one can choose a sequence $\{k_n\}_{n=1}^\infty$ such that $k_n\to\infty$ as $n\to\infty$ and the supports of the functions $\psi_{k_n}$ are mutually disjoint which means that $e+ p^2\in \sigma_\mathrm{ess}(H)$ for all $p\in \mathbb{R}$, and, one infers that $[e, \infty) \subset \sigma_{\mathrm{ess}}(H)$.

Next we must establish that the spectrum of $H$ below $e$, if any, can only be discrete. The Neumann bracketing gives the estimate
\begin{equation}
\label{bracketing}
H\ge H_1\oplus H_2,
\end{equation}
where $H_1$ is the Neumann restriction of $H$ on $L^2(\mathbb{R}\times \{|x|> c\})$, i.e., the operator with the Neumann boundary conditions at $x=\pm c$, and $H_2$ is the complementary Neumann restriction on $L^2(\{|x|< c\})$. We will show that both operators $H_1$ and $H_2$ have purely discrete spectrum below $e$. Then the same will hold for the direct sum $H_1 \oplus H_2$ and via the minimax principle and (\ref{bracketing}) also for the operator $H$.

Now, just by the fact that $V(x, y)= \tilde{V}(y)$ and $A(x, y)=0$ if $|x|>c$, one can easily check that
\begin{align*}
 &\int_{\mathbb{R}\times \{|x|>c\}} |i \nabla u + A u|^2\,d x\,d y-  \int_{\mathbb{R}\times \{|x|>c\}} V(x, y)|u|^2\,d x\,d y \\ & \ge
\int_{\{|x|>c\}} \int_{\mathbb{R}}|u_y|^2\,d y\,d x+ \int_{\{|x|>c\}} \int_{\mathbb{R}}|u_x|^2\,d x\,d y-  \int_{\{|x|>c\}} \int_{\mathbb{R}} \tilde{V}(y) |u|^2\,d x\,d y 
\\ &
\ge\int_{\{|x|>c\}} \int_{\mathbb{R}}|u_y|^2\,d yd x-  \int_{\{|x|>c\}} \int_{\mathbb{R}} \tilde{V}(y) |u|^2\,d x\,d y.\\ 
\end{align*}

Since the principal eigenvalue of the operator $-\frac{\mathrm{d}^2}{\mathrm{d} t^2}- \tilde{V}$ is $e$, the above estimate implies
\begin{align*}
\int_{\mathbb{R}\times \{|x|>c\}} |i \nabla u + A u|^2\,d xd y-  \int_{\mathbb{R}\times \{|x|>c\}} V(x, y)|u|^2\,d xd y&\\\ge e \int_{\mathbb{R}\times \{|x|>c\}}|u|^2\,d xd y\,;
\end{align*}
this means that the spectrum of $H_1$ below $e$ is empty.

It remains to deal with the operator $H_2$. Again using  the Neumann bracketing method, we get 
\begin{equation}\label{simple}H_2\ge H^2_1 \oplus H^2_2,\end{equation}
where $H^2_1$ is the Neumann restriction of $H_2$ to $L^2((-c, c)\times (0, \alpha))$ with
$\alpha:= \underset{x\in (-c, c)}{\mathrm{max}}\,|V(x, y)|$ and $H^2_2$ is the Neumann restriction of $H_2$ to $L^2((-c, c)\times \{|y|>\alpha\})$.

One has the simple estimate 
$$H^2_1\ge (i \nabla+ A)^2- \alpha,
$$
which together with the discreteness of the magnetic Neumann Laplacian on the rectangle $(-c, c)\times (-\alpha, \alpha)$
implies the discreteness of the spectrum of $H^2_1$.

Moving on to the operator $H_2^2$, one can easily see that it is non-negative due to the fact that on its domain potential $V$ is zero. Thus the negative spectrum of $H_2^2$ is absent.

Therefore, in view of (\ref{simple}), we establish that the negative spectrum of $H_2$, if it exists, consists of the finite number of eigenvalues of finite multiplicity. This together with (\ref{bracketing}) and the fact that the spectrum of $H_1$ below $e$ is empty, completes the proof of our claim.

\bigskip

Let us now pass to the proof of Theorem~\ref{th.22}. We need the following lemma (the proof is given in the Appendix):
\begin{lemma}\label{Appendix} Let $L:= -\frac{d^2}{d t^2}+Q$, where $Q$ is a compactly supported function belonging to $L^\infty(\mathbb{R})$. 
Let $\lambda_1$ be the ground state eigenvalue of $L$. Then the corresponding eigenfunction is nowhere
zero (it can be chosen to be positive).
\end{lemma}

Let $g$ be the ground state eigenfunction of the operator $h$. Given the above lemma, it can be chosen to be positive. Thus we represent each function $\psi\in C_0^\infty(\mathbb{R}^2)$ as a multiplication $\varphi(x, y) g(y)$ with $\varphi\in C_0^\infty(\mathbb{R}^2)$. Then with the gauge 
$A(x, y)=\Big(-\int_0^y B(x, t)\,\mathrm{d}t,0\Big)$, for any $\psi\in C_0^\infty(\mathbb{R}^2)$  we have
\begin{eqnarray}\nonumber&
\int_{\mathbb{R}^2}\left(|i \nabla \psi(x, y)+ A(x, y) \psi(x, y)|^2- V(x, y) |\psi(x, y)|^2\right)\,d x\,d y\\\nonumber&=
\int_{\mathbb{R}^2}\left|i \frac{\partial \varphi}{\partial x}(x, y)- \varphi(x, y) \int_0^x B(t, y)\,d t\right|^2 g(y)^2\,d x\,d y+ \int_{\mathbb{R}^2}\left|\frac{\partial \varphi}{\partial y}(x, y) g(y)+  \varphi(x, y) g'(y)\right|^2\,d x\,d y\\\nonumber&- \int_{\mathbb{R}^2}
V(x, y) |\varphi(x, y)|^2 g(y)^2\,d x\,d y\\\nonumber& =\int_{\mathbb{R}^2}\left|i \frac{\partial \varphi}{\partial x}(x, y)- \varphi(x, y) \int_0^x B(t, y)\,d t\right|^2 g(y)^2\,d x\,d y\\\nonumber&+  \int_{\mathbb{R}^2}\left(\left|\frac{\partial \varphi}{\partial y}(x, y)\right|^2 g(y)^2+ \frac{\partial \varphi}{\partial y}(x, y) \overline{\varphi}(x, y) g(y) g'(y)+ \varphi(x, y) \frac{\partial \overline{\varphi}}{\partial y}(x, y) g(y) g'(y)+|\varphi(x, y)|^2 g'(y)^2\right)\,d xdy \\\label{calc.}&- \int_{\mathbb{R}^2}V(x, y) |\varphi(x, y)|^2 g(y)^2\,d x\,d y.
\end{eqnarray}

One can easily check by integrating by parts that 
\begin{eqnarray*}
\int_{\mathbb{R}^2}\left(\frac{\partial \varphi}{\partial y}(x, y) \overline{\varphi}(x, y) g(y) g'(y)+ \varphi(x, y) \frac{\partial \overline{\varphi}}{\partial y}(x, y) g(y) g'(y)+|\varphi(x, y)|^2 g(y)'^2\right)\,d xd y\\= -\int_{\mathbb{R}^2}g''(y)g(y) |\varphi(x, y)|^2\,d xd y\,,
\end{eqnarray*}
which together with (\ref{calc.}) and the fact that $g$ is the ground state eigenfunction of $h$ gives
\begin{align}\nonumber&
\int_{\mathbb{R}^2}\left(|i \nabla \psi(x, y)+ A(x, y) \psi(x, y)|^2- V(x, y) |\psi(x, y)|^2\right)\,d x\,d y\\\nonumber& =\int_{\mathbb{R}^2}\left|i \frac{\partial \varphi}{\partial x}(x, y)- \varphi(x, y) \int_0^x B(t, y)\,d t\right|^2 g(y)^2\,d x\,d y+  \int_{\mathbb{R}^2}\left|\frac{\partial \varphi}{\partial y}(x, y)\right|^2 g(y)^2\,d x\,d y\\\nonumber& +\int_{\mathbb{R}^2}(-g''(y) -V(x, y)g(y)) g(y) |\varphi(x, y)|^2\,d x\,d y\\\nonumber&= \int_{\mathbb{R}^2}\left|i \frac{\partial \varphi}{\partial x}(x, y)- \varphi(x, y) \int_0^x B(t, y)\,d t\right|^2 g(y)^2\,d x\,d y+  \int_{\mathbb{R}^2}\left|\frac{\partial \varphi}{\partial y}(x, y)\right|^2 g(y)^2\,d x\,d y\\\nonumber& +\int_{\mathbb{R}^2}(-g''(y) -\tilde{V}(y) g(y)) g(y) |\varphi(x, y)|^2\,d xd y+ \int_{\mathbb{R}^2}(\tilde{V}(y)-V(x, y)) g(y)^2 |\varphi(x, y)|^2\,d x\,d y \\\nonumber&= \int_{\mathbb{R}^2}\left|i \frac{\partial \varphi}{\partial x}(x, y)- \varphi(x, y) \int_0^x B(t, y)\,d t\right|^2 g(y)^2\,d x\,d y+  \int_{\mathbb{R}^2}\left|\frac{\partial \varphi}{\partial y}(x, y)\right|^2 g(y)^2\,d x\,d y\\\nonumber&+ \int_{\mathbb{R}^2}(\tilde{V}(y)-V(x, y)) g(y)^2 |\varphi(x, y)|^2\,d x\,d y +e \int_{\mathbb{R}^2}g(y)^2 |\varphi(x, y)|^2 \,d x\,d y\\\nonumber& \ge\int_{\mathbb{R}^2}\left|i \frac{\partial \varphi}{\partial x}(x, y)- \varphi(x, y) \int_0^x B(t, y)\,d t\right|^2 g(y)^2\,d x\,d y+  \int_{\mathbb{R}^2}\left|\frac{\partial \varphi}{\partial y}(x, y)\right|^2 g(y)^2\,d x\,d y\\\label{expression}&+ \int_{\Omega^a}(\tilde{V}(y)-V(x, y)) g(y)^2 |\varphi(x, y)|^2\,d x\,d y +e \int_{\mathbb{R}^2}g(y)^2 |\varphi(x, y)|^2 \,d x\,d y\,.
\end{align}

We introduce a new constant $\tau>0$ so that the ball $\mathbb{B}(0, \tau)$ with center at zero and radius $\tau$ contains $\Omega^a$.
Let $\beta_g:=\underset{|z|\le \tau}{\mathrm{min}}\,g(z)^2$. From the positivity of $g$ it follows that $\beta_g\ne0$. Then
estimate (\ref{expression}) implies
\begin{eqnarray}\nonumber
\int_{\mathbb{R}^2}\left(|i \nabla \psi(x, y)+ A(x, y) \psi(x, y)|^2- V(x, y) |\psi(x, y)|^2\right)\,d x\,d y\\\nonumber\ge \beta_g\left(\int_{\mathbb{B}(0, \tau)}\left|i \frac{\partial \varphi}{\partial x}(x, y)- \varphi(x, y) \int_0^x B(t, y)\,d t\right|^2\,d x\,d y+  \int_{\mathbb{B}(0, \tau)}\left|\frac{\partial \varphi}{\partial y}(x, y)\right|^2\,d x\,d y\right)\\\label{Hardy}- \beta_g \|\tilde{V}(y)-V(x, y)\|_{L^\infty(\Omega^a)}\,\int_{\Omega^a} |\varphi(x, y)|^2\,d x\,d y +e \int_{\mathbb{R}^2} g(y)^2 |\varphi(x, y)|^2\,d x\,d y\,.
\end{eqnarray}

To continue the proof, we need the following lemma:
\begin{lemma}\label{2lemma}
Let $I\subset \mathbb{R}$ be an interval and denote its length by $|I|$, then for any $g\in \mathcal{H}^1(I)$ and all $x\in I$ we have
$$
|g(x)|^2\le 2|I| \int_I |g'(t)|^2\,d t +\frac{2}{|I|} \int_I |g(t)|^2\,d t
$$
\end{lemma}
\begin{proof}
It is easy to check that there exists a point $x_0\in I$ such that
$$
|g(x_0)|\le \frac{1}{\sqrt{|I|}} \sqrt{\int_I |g(t)|^2\,d t}\,;
$$
this fact, combined with the inequalities of Schwarz and Jensen, yields
$$|g(x)|^2=\left|\int_{x_0}^x g'(t)\,d t- g(x_0)\right|^2\le 2|I| \int_I |g'(t)|^2\,d t +\frac{2}{|I|} \int_I |g(t)|^2\,d t,$$
which is what we have set up to prove.
\end{proof}
 
Using the fact that for any fixed $x\in(-c, c)$ function $\varphi(x, y)$ belongs to the Sovolev space 
$$\mathcal{H}^1\left(-\sqrt{\tau^2-x^2}, \sqrt{\tau^2-x^2}\right)$$ and by Lemma(\ref{2lemma}) we conclude that
\begin{equation}\label{above}
|\varphi(x, y)|^2\le 4\sqrt{\tau^2-x^2} \int_{-\sqrt{\tau^2-x^2}}^{\sqrt{\tau^2-x^2}} \left|\frac{\partial{\varphi}}{{\partial y}}(x, t)\right|^2\,d t
+\frac{1}{\sqrt{\tau^2-x^2}} \int_{-\sqrt{\tau^2-x^2}}^{\sqrt{\tau^2-x^2}} |\varphi(x, t)|^2\,d t\,.
\end{equation}

Let 
$$
\Omega^a\cap\mathcal{B}(0, c)=\left\{x\in(-c, c), \gamma_1(x)<y<\gamma_2(x)\right\},
$$
where $\gamma_1:(-c, c)\to\mathbb{R}$ and $\gamma_2:(-c, c)\to\mathbb{R}$ are smooth curves, be the representation of the set $\Omega^a\cap\mathcal{B}(0, c)$. Then
\begin{eqnarray*}
\int_{\Omega^a\cap\mathcal{B}(0, c)}|\varphi(x, y)|^2\,d x\,d y=
\int_{-c}^c \int_{\gamma_1(x)}^{\gamma_2(x)}|\varphi(x, y)|^2\,d x\,d y\\
=\int_{-c}^c \int_{\gamma_1(x)}^{\gamma_2(x)} \biggr( 4\sqrt{\tau^2-x^2} \int_{-\sqrt{\tau^2-x^2}}^{\sqrt{\tau^2-x^2}} \left|\frac{\partial{\varphi}}{{\partial y}}(x, t)\right|^2\,d t +\frac{1}{\sqrt{\tau^2-x^2}}\int_{-\sqrt{\tau^2-x^2}}^{\sqrt{\tau^2-x^2}} |\varphi(x, t)|^2\,d t\biggl)\,d x\,d y
\\\le\int_{-c}^c \int_{\gamma_1(x)}^{\gamma_2(x)} \biggr( 4\tau \int_{-\sqrt{\tau^2-x^2}}^{\sqrt{\tau^2-x^2}} \left|\frac{\partial{\varphi}}{{\partial y}}(x, t)\right|^2\,d t +\frac{1}{\sqrt{\tau^2-c^2}} \int_{-\sqrt{\tau^2-x^2}}^{\sqrt{\tau^2-x^2}} |\varphi(x, t)|^2\,d t\biggl)\,d x\,d y\\\le 4\tau \|\gamma_1-\gamma_2\|_{L^\infty(-c, c)} \int_{{\mathbb{B}(0, \tau)}}\left|\frac{\partial{\varphi}}{{\partial y}}(x, t)\right|^2\,d t\,d x
+\frac{ \|\gamma_1-\gamma_2\|_{L^\infty(-c, c)}}{\sqrt{\tau^2-c^2}} \int_{\mathbb{B}(0, \tau)}|\varphi(x, t)|^2\,d t\,d x\\
\le 4\tau \delta(\Omega^a) \int_{{\mathbb{B}(0, \tau)}}\left|\frac{\partial{\varphi}}{{\partial y}}(x, t)\right|^2\,d t\,d x
+\frac{\delta(\Omega^a)}{\sqrt{\tau^2-c^2}} \int_{\mathbb{B}(0, \tau)}|\varphi(x, t)|^2\,d t\,d x\,.
\end{eqnarray*}

Applying the above inequality to the right-hand side of (\ref{Hardy}), we get
\begin{eqnarray}\nonumber
\int_{\mathbb{R}^2}\left(|i \nabla \psi(x, y)+ A(x, y) \psi(x, y)|^2- V(x, y) |\psi(x, y)|^2\right)\,d xd y\\\nonumber\ge \beta_g\left(\int_{\mathbb{B}(0, \tau)}\left|i \frac{\partial \varphi}{\partial x}(x, y)- \varphi(x, y) \int_0^x
B(t, y)\,d t\right|^2\,d x\,d y+  \int_{\mathbb{B}(0, \tau)}\left|\frac{\partial \varphi}{\partial y}(x, y)\right|^2\,d x\,d y\right)\\\nonumber- \beta_g\,\|\tilde{V}(y)-V(x, y)\|_{L^\infty(\Omega^a)}\,\left( 4\tau\delta(\Omega^a) \int_{{\mathbb{B}(0, \tau)}}\left|\frac{\partial{\varphi}}{{\partial y}}(x, y)\right|^2\,d x\,d y+ \frac{\delta(\Omega^a)}{\sqrt{\tau^2-c^2}} \int_{\mathbb{B}(0, \tau)}|\varphi(x, y)|^2\,d x\,d y\right)\\\label{Hardy1}+e \int_{\mathbb{R}^2} g(y)^2 |\varphi(x, y)|^2\,d x\,d y\,.
\end{eqnarray}

Under the assumption that 
\begin{equation}\label{1assump.}
\|\tilde{V}(y)-V(x, y)\|_{L^\infty(\Omega^a)}< \frac{1}{8\tau\delta(\Omega^a) }\end{equation}
inequality (\ref{Hardy1}) implies
\begin{eqnarray}\nonumber
\int_{\mathbb{R}^2}\left(|i \nabla \psi(x, y)+ A(x, y) \psi(x, y)|^2- V(x, y) |\psi(x, y)|^2\right)\,d x\,d y\\\nonumber\ge \beta_g
\int_{\mathbb{B}(0, \tau)}\left|i \frac{\partial \varphi}{\partial x}(x, y)- \varphi(x, y) \int_0^x B(t, y)\,d t\right|^2\,d x\,d y+  \frac{\beta_g}{2}
\int_{\mathbb{B}(0, \tau)}\left|\frac{\partial \varphi}{\partial y}(x, y)\right|^2\,d x\,d y\\\nonumber- \frac{\delta(\Omega^a)\,\beta_g}{\sqrt{\tau^2-c^2}} \|\tilde{V}(y)-V(x, y)\|_{L^\infty(\Omega^a)}\,\int_{\mathbb{B}(0, \tau)}|\varphi(x, y)|^2\,d x\,d y+e \int_{\mathbb{R}^2} g(y)^2 |\varphi(x, y)|^2\,d x\,d y\\\nonumber\ge \frac{\beta_g}{2}
 \int_{\mathbb{B}(0, \tau)}\left(|i \nabla \varphi(x, y)+ A(x, y) \varphi(x, y)|^2\right)\,d x\,d y
\\\nonumber- \frac{\delta(\Omega^a)\,\beta_g}{\sqrt{\tau^2-c^2}} \|\tilde{V}(y)-V(x, y)\|_{L^\infty(\Omega^a)}\int_{\mathbb{B}(0, \tau)}|\varphi(x, y)|^2\,d x\,d y\\\label{Hardy2}+e \int_{\mathbb{R}^2} g(y)^2 |\varphi(x, y)|^2\,d x\,d y\,.
\end{eqnarray}

Using the Hardy inequality for magnetic Dirichlet forms \cite{EK052}
$$
\int_{\mathbb{B}(0, \tau)}|i \nabla \varphi+A \varphi|^2\,d xd y\ge \mathcal{M} \int_{\mathbb{B}(0, \tau)}|\varphi|^2\,d x\,d y,
$$
where the constant $\mathcal{M}>0$ depends on the magnetic potential $B$, the right-hand side of (\ref{Hardy2}) can be estimated as follows

\begin{align}\nonumber&
\int_{\mathbb{R}^2}\left(|i \nabla \psi(x, y)+ A(x, y) \psi(x, y)|^2- V(x, y) |\psi(x, y)|^2\right)\,d x\,d y&\\\nonumber&\ge \left(\frac{\beta_g}{2}
\mathcal{M}- \frac{\delta(\Omega^a)\,\beta_g}{\sqrt{\tau^2-c^2}} \|\tilde{V}(y)-V(x, y)\|_{L^\infty(\Omega^a)}\right)  \int_{\mathbb{B}(0, \tau)}|\varphi(x, y)|^2\,d xd y&\\\nonumber&+ e \int_{\mathbb{R}^2} g(y)^2 |\varphi(x, y)|^2\,d x\,d y\,.&\end{align}

Thus, one can easily check that the assumption
$$
 \|\tilde{V}(y)-V(x, y)\|_{L^\infty(\Omega^a)}< \frac{
\mathcal{M} \sqrt{\tau^2-c^2}}{2\delta(\Omega^a)}
$$
guarantees that the right-hand side of the above inequality is not less than $e \int_{\mathbb{R}^2} g(y)^2 |\varphi(x, y)|^2\,
d x\,d y$.

Since we also used the assumption (\ref{1assump.}) before, the statement of the theorem is established with
$$
\|\tilde{V}(y)-V(x, y)\|_{L^\infty(\Omega^a)}\le\mathcal{C} \left(\delta(\Omega^a)\right)^{-1},
$$
where
$$\mathcal{C}= \mathrm{min}\left\{ \frac{
 \mathcal{M} \sqrt{\tau^2-c^2}}{2},\,\frac{1}{8\tau}\right\}.$$
 
On the other hand, using the fact that $V(x, y)=\tilde{V}(u(x, y))$, where $u$ is the distance from $(x, y)$ to $\Gamma$ as usual, we complete the proof of the theorem.  

\bigskip
The next our remark presents a concrete example when the magnetic field destroys the eigenvalues resulting from the curve geometry.
  
\begin{remark}\label{final}

In the paper \cite{E20} it was proved that if $\Gamma$ is a $C^2$- smooth curve, which is a local deformation of a straight line $\{(x, 0)\}_{x\in\mathbb{R}}$, then the operator $-\Delta- \frac{1}{\varepsilon}V\left(\frac{u}{\varepsilon}\right)$, where $u$ is the distance from a point to $\Gamma$ and a nonzero potential $V\ge0$ from $L^\infty(\mathbb{R})$ with $\mathrm{supp}V\subset[-a, a]$, has a non-empty discrete spectrum provided that $\varepsilon$ is small enough.
 
In terms of our notations $\tilde{V}(z)= \frac{1}{\varepsilon} V\left(\frac{z}{\varepsilon}\right),\,z\in(-\varepsilon a, \varepsilon a)$.

One can choose a suitable smooth curve $\Gamma$ in order to have the asymptotics $d(\Omega^\varepsilon)= \mathcal{O}(\varepsilon)$.
Then if we choose a potential $V$ being sufficiently small in order to guarantee the validity of (\ref{assumption1}), and add a suitable magnetic field, the assumptions of Theorem \ref{th.11.} are easily satisfied
and therefore the discrete spectrum of the operator $(i \nabla+A)^2- \frac{1}{\varepsilon}V\left(\frac{u}{\varepsilon}\right)$ becomes empty.

\end{remark}

\section{Appendix}

In the proof of Lemma \ref{Appendix} we use the methods of (\cite{K19}) and (\cite{K10}).
If $v$ is a normalized eigenfunction of $L$ corresponding to the ground state eigenvalue $\lambda_1$, then it follows from the Rayleigh-Ritz formula \cite{RS78}
$$
\lambda_1 =\int_{\mathbb{R}}|\nabla v|^2\, d t+\int_{\mathbb{R}}Q |v|^2\,d t
$$
that $|v|\ge 0$ is also one. We will now prove that it is  never equal to zero. This would imply that the eigenfunction can be chosen to be strictly positive. We have to exclude that $v$ can have zeroes. Suppose that $x=x_0$ is the zero of $v$.
Since the function $v$ is continuously differentiable in particular when $x = x_0$,
then its derivative there should be equal to zero, since the function is non-negative
and $v(x_0) = 0$. It follows that at this particular point the function $v$ satisfies zero Cauchy data
$$ \begin{cases}v(x_0) = 0\\ 
 v'(x_0)=0\end{cases}
$$
and is therefore identically zero on the whole $\mathbb{R}$ as a solution to the
ordinary second-order differential equation. This implies that $v$ should be identically zero. This concludes our proof.

\subsection*{Acknowledgements}

J.B.R. acknowledges support from Instituto Politecnico Nacional, grant number SIP20230312. D.S. acknowledges the support of the Czech Science Foundation (GACR), project 21-07129S. D.S. and B.S. acknowledge the support from the Czech-Polish project \newline BPI/PST/2021/1/00031.


\begin{thebibliography}{10}

\bibitem{BE21} D.~ Barseghyan, P.~ Exner,
Magnetic field influence on the discrete spectrum of locally deformed leaky wires,
Rep. Math. Phys. 88 (2021), 47--57.
 
\bibitem{EK05} T.~Ekholm, H.~Kova\v{r}\'{\i}k, Stability of the magnetic Schr\"{o}dinger operator in a waveguide, Comm. PDE  30 (2005), 539--565.

\bibitem{EK15} P.~Exner, H.~Kova\v{r}\'{\i}k: Quantum Waveguides, Springer International, Heidelberg 2015.

\bibitem{E20} P.~Exner, Spectral properties of soft quantum waveguides, 
Journal of Physics A: Mathematical and Theoretical. 53: 35 (2020).
 
\bibitem{EKK08}
 T.~Ekholm, H.~Kova\v{r}\'{\i}k, D.~Krej\v{c}i\v{r}\'{\i}k,
A Hardy inequality in twisted waveguides,
 Arch. Rational Mech. Anal. 188 (2008), 245--264.

\bibitem{EK051}
 P.~Exner, H.~Kova\v{r}\'{\i}k, Spectrum of the Schr\"{o}dinger operator in a perturbed periodically twisted tube,
Lett. Math. Phys. 73 (2005), 183--192.

\bibitem{EK052}
T.~Ekholm, H.~Kova\v{r}\'{\i}k, Stability of the Magnetic Schr\"{o}dinger Operator in a Waveguide, Communications in Partial Differential Equations, 30: 539–565, 2005.

\bibitem{ES89}
P.~Exner, P.~\v{S}eba,
Bound states in curved quantum wavequides, J. Math. Phys. 30 (1989), 2574--2580.

\bibitem{K10} D.~Krej\v{c}i\v{r}\'{\i}k, Schr\"{o}dinger operators and their spectra, lecture notes (2010).

\bibitem{K19} P.~Kurasov, On the ground state for quantum graphs, Letters in Mathematical Physics (2019) 109: 2491--2512.

\bibitem{RS81} M.~Reed, B.~Simon, Methods of Modern Mathematical Physics, I. Functional Analysis, Academic Press, New York 1981.

\bibitem{RS78} M.~Reed, B.~Simon, Methods of modern mathematical physics, IV. Analysis of opera-
tors, Academic Press, New York, 1978.

\end{thebibliography}
\end{document}